\documentclass[12pt]{article}
\usepackage{amssymb}
\usepackage{amsmath}
\usepackage{amsfonts}
\usepackage{mathrsfs}

\newcommand\qed{\hfill $\square$}

\newcommand{\NN}{\mathbb N}
\newcommand{\RR}{\mathbb R}

\newcommand{\CC}{\mathbb C}

\newcommand{\dd}{\,\mathrm{d}}

\newtheorem{theo}{Theorem}
\newtheorem{prop}[theo]{Proposition}
\newtheorem{lem}[theo]{Lemma}
\newtheorem{cor}[theo]{Corollary}
\newtheorem{rem}[theo]{Remark}

\newcommand{\beqn}{\begin{equation}}
\newcommand{\eeqn}{\end{equation}}
\newcommand{\bear}{\begin{eqnarray}}
\newcommand{\eear}{\end{eqnarray}}
\newcommand{\bean}{\begin{eqnarray*}}
\newcommand{\eean}{\end{eqnarray*}}

\begin{document}
\begin{center}


{\LARGE Asymptotic behavior for a viscous\\[1ex]
Hamilton-Jacobi equation with critical exponent}
\\[12mm]
{\bf Thierry Gallay}\\[1mm]
Institut Fourier\\
Universit{\'e} de Grenoble I\\
B.P. 74 \\
38402 Saint-Martin-d'H{\`e}res, France\\
{\tt Thierry.Gallay@ujf-grenoble.fr}
\\[8mm]
{\bf Philippe Lauren\c{c}ot}\\[1mm]
Math\'ematiques pour l'Industrie et la Physique, CNRS UMR~5640\\ 
Universit\'e Paul Sabatier -- Toulouse~3\\
118, route de Narbonne\\
31062 Toulouse cedex 9, France\\
{\tt laurenco@mip.ups-tlse.fr}
\end{center}



\vspace*{1cm}

\begin{abstract}
The large time behavior of non-negative solutions to the viscous
Hamilton-Jacobi equation $\partial_t u - \Delta u + |\nabla
u|^q = 0$ in $(0,\infty)\times\RR^N$ is investigated for the
critical exponent $q=(N+2)/(N+1)$. Convergence towards a rescaled
self-similar solution to the linear heat equation is shown, the
rescaling factor being $(\ln{t})^{-(N+1)}$. The proof relies on the
construction of a one-dimensional invariant manifold for a suitable
truncation of the equation written in self-similar variables.
\end{abstract}

\thispagestyle{empty}

\vfill 


\newpage

\section{Introduction} \label{int}

The dynamics of integrable non-negative solutions to the viscous 
Hamilton-Jacobi equation 
\beqn
\label{vhj}
\partial_t u - \Delta u + |\nabla u|^q \,=\, 0\,,
\quad (t,x)\in (0,\infty)\times\RR^N\,,
\eeqn
depends strongly on the value of the parameter $q\in (0,\infty)$ and
results from the competition between the linear diffusion term $\Delta
u$ and the nonlinear absorption term $|\nabla u|^q$. An
important issue is therefore to determine which mechanism (diffusion
or absorption) is dominant for large times. A first indication is
given by the behavior of the $L^1$ norm $\|u(t)\|_{L^1}$,
which is time-independent for non-negative solutions of the heat
equation and strictly decreasing for nontrivial non-negative
solutions of (\ref{vhj}). For such solutions, it is proved in
\cite{AB98,BL99,BK99} that
\beqn
\label{liml1}
I_\infty \,:=\, \lim_{t\to\infty} \|u(t)\|_{L^1} 
\left\{
\begin{array}{lcl}
> 0 & \mbox{ if } & q>q_\star:=\displaystyle{\frac{N+2}{N+1}}\,, \\
 & & \\
= 0 & \mbox{ if } & q\in (0,q_\star]\,. 
\end{array}
\right.
\eeqn
This suggests that diffusion dominates the large time behavior when
$q>q_\star$, whereas absorption becomes effective for $q\le q_\star$. 
As a matter of fact, if $q>q_\star$ it is shown in
\cite{BKaL04,BGK04} that the nonlinear term $|\nabla u|^q$ 
becomes negligible for large times, and that the solution of 
(\ref{vhj}) behaves as $t \to \infty$ like the self-similar solution
$I_\infty\,g$ of the linear heat equation, where
\beqn
\label{gk}
g(t,x) \,=\, \frac{1}{t^{N/2}}\,G\left( \frac{x}{t^{1/2}} \right)
\;\;\mbox{ and }\;\; G(\xi) \,=\, \frac{1}{(4\pi)^{N/2}}\,\exp{\left( -
\frac{|\xi|^2}{4} \right)}\,. 
\eeqn
On the other hand, if $q\in (1,q_\star)$, both diffusion and 
absorption play a role in the large time asymptotics. Indeed, 
if $u(0,x)$ decays faster than $|x|^{-\alpha}$ as $|x|\to \infty$ 
with $\alpha = (2-q)/(q-1) > N$, it is proved in \cite{BKaL04}
that the solution $u(t)$ converges as $t \to \infty$ to the 
so-called \textit{very singular solution}, a self-similar solution
of (\ref{vhj}) whose existence and uniqueness have been established 
in \cite{BKoL04,BL01,QW01}. In that case, the $L^1$-norm of $u(t)$ 
decays to zero like $t^{-(\alpha-N)/2}$ as $t \to \infty$. Finally, 
the influence of the absorption term $|\nabla u|^q$ is 
much stronger for $q\in (0,1]$: depending on the initial data, 
one might have exponential decay of the solution as $t \to \infty$ 
\cite{AB98,BRV96,BRV97,Gi05}, or even extinction in finite time 
if $q\in (0,1)$ \cite{BLS02,BLSS02,Gi05}. For such values of the 
parameter, it is the diffusion term which is expected to be
negligible for large times. 

To summarize, precise asymptotic expansions show that the large time
behavior of non-negative solutions to (\ref{vhj}) with sufficiently
localized initial data is determined by the sole diffusion if
$q>q_\star$, whereas absorption plays an important role if
$q<q_\star$. With this perspective in mind, it is interesting to
investigate the critical case $q=q_\star=(N+2)/(N+1)$ where a
transition between both regimes is expected to occur. Very few results
are available in this situation: we only know that $\|
u(t)\|_{L^1}\to 0$ as $t\to\infty$, as already stated in
(\ref{liml1}), and that $\|u(t)\|_{L^1}$ cannot decay faster
than $(\ln{t})^{-(N+1)}$ for large times \cite[Proposition~3]{BLS02}.
The purpose of this work is to fill this gap and to give an
accurate description of the large time behavior of the non-negative
solutions to
\bear
\label{vhj1}
\partial_t u - \Delta u + |\nabla u|^{q_\star} & = & 0\,,
\quad (t,x)\in (0,\infty)\times\RR^N\,,\\
\label{vhj2}
u(0) & = & u_0\,, \quad x\in\RR^N\,,
\eear
when the initial data $u_0(x)$ decay to zero sufficiently rapidly 
as $|x|\to\infty$. More precisely, we assume that 
$u_0 \ge 0$ belongs to the weighted $L^2$ space
\beqn
\label{L2mdef}
L^2_m(\RR^N) \,=\, \Bigl\{u \in L^2(\RR^N) \,\Big|\, 
|u|_m := \|(1{+}|x|^{2m})^{1/2}u\|_{L^2} < \infty\Bigr\}\,,
\eeqn
for some $m > N/2$. Then (by H\"older's inequality) $u_0 \in
L^1(\RR^N) \cap L^2(\RR^N)$ and it follows from \cite{BL99,BSW02,GGK03}
that the Cauchy problem (\ref{vhj1}), (\ref{vhj2}) has a unique 
global solution $u \in \mathcal{C}([0,\infty);L^1(\RR^N)) \cap
\mathcal{C}((0,\infty);W^{1,\infty}(\RR^N))$. Our main result
describes the large time behavior of this solution: 

\begin{theo}\label{th1}
Assume that the initial condition $u_0$ is non-negative, not identically 
zero, and belongs to $L^2_m(\RR^N)$ for some $m > N/2$. Then the
(unique) solution $u$ to (\ref{vhj1}),(\ref{vhj2}) satisfies, 
for all $p\in [1,\infty]$, 
\beqn
\label{dr}
\lim_{t\to\infty} t^{\frac{N}2(1-\frac1p)}\,(\ln{t})^{N+1}\,\left\|u(t) -
\frac{M_\star}{(\ln{t})^{N+1}}\,g(t) \right\|_{L^p} = ~0\,, 
\eeqn 
where $M_\star = (N+1)^{N+1}\,\|\nabla G\|_{L^{q_\star}}^{-(N+2)}$
and $g(t,x)$, $G(\xi)$ are defined in (\ref{gk}). 
\end{theo}

In other words, if the initial condition decays sufficiently rapidly at 
infinity, the solution $u$ to (\ref{vhj1}) behaves asymptotically 
like a particular self-similar solution $M_\star\,g$ of the linear 
heat equation, with an extra logarithmic factor due to the effect 
of the absorption term. Such a logarithmic correction also appears
in other parabolic equations with absorption and critical exponent,
for instance in the nonlinear diffusion equation $\partial_t u - 
\Delta u^m + u^{m+(2/N)}=0$ with $m\ge 1$, see
e.g. \cite{GH04,GV91,He99} and the references therein. 

\begin{rem}\label{re1}
As $\|g(t)\|_{L^1}=1$ for all $t>0$, the $L^1$-norm of the
solution $u(t)$ behaves exactly like $M_\star(\ln{t})^{-(N+1)}$ for
large times under the assumptions of Theorem~\ref{th1}. This has to
be compared with \cite[Proposition~3]{BLS02}, where it is shown that
there is no nontrivial non-negative solution of (\ref{vhj1}) such
that $\|u(t)\|_{L^1} \le C(\ln{t})^{-\gamma}$ for $\gamma > N+1$.
In fact, using Theorem~\ref{th1} and a comparison argument, it is
straightforward to verify that, for all nontrivial non-negative
integrable data, the solution of (\ref{vhj1}) satisfies
$$
\liminf_{t\to\infty}\,(\ln{t})^{N+1}\|u(t)\|_{L^1} \,\ge\, M_\star\,.
$$
\end{rem}

\medskip

Our analysis of the large time behavior of solutions to (\ref{vhj1}),
(\ref{vhj2}) relies on an alternative formulation of (\ref{vhj1}) in
terms of the so-called ``scaling variables'' or ``similarity
variables''  
\beqn
\label{scavar}
\xi \,=\, \frac{x}{(1+t)^{1/2}} \;\;\mbox{ and }\;\; \tau  \,=\, 
\ln{(1+t)}\,.
\eeqn
Introducing the new unknown function $v$ defined by
\beqn
\label{newuf}
u(t,x) \,=\, \frac{1}{(1+t)^{N/2}}\,v\left( \ln{(1+t)} ,
\frac{x}{(1+t)^{1/2}} \right)\,, \quad (t,x)\in
[0,\infty)\times\RR^N\,, 
\eeqn
we deduce from (\ref{vhj1}), (\ref{vhj2}) that $v(\tau,\xi)$ solves 
the initial-value problem 
\bear
\label{rsvhj1}
\partial_\tau v & = & \mathcal{L} v - |\nabla v|^{q_\star}\,,
\quad (\tau,\xi)\in (0,\infty)\times\RR^N\,, \\ 
\label{rsvhj2}
v(0) & = & u_0\,, \quad \xi\in\RR^N\,,
\eear
where the linear operator $\mathcal{L}$ is given by
\beqn
\label{opL}
\mathcal{L} v(\xi) \,=\, \Delta v(\xi) + \frac{1}{2}\, \xi\cdot\nabla
v(\xi) + \frac{N}{2}\, v(\xi)\,, \quad \xi\in\RR^N\,. 
\eeqn
Observe that equation (\ref{rsvhj1}) is still autonomous, although
it was obtained from (\ref{vhj1}) through the time-dependent 
transformation (\ref{newuf}). This crucial property follows 
from the fact that (\ref{vhj1}) is invariant under the rescaling
$u(t,x) \mapsto \lambda^N u(\lambda^2 t,\lambda x)$, because 
$q_\star=(N+2)/(N+1)$. Remark also that $\mathcal{L} G = 0$, where 
$G$ is defined in (\ref{gk}). 

At this stage, we follow the approach of \cite{GW02,Wa97} and prove
that the large time behavior of the solutions of (\ref{rsvhj1}),
(\ref{rsvhj2}) is governed, up to exponentially decaying terms, 
by an ordinary differential equation which results from restricting 
the dynamics of (\ref{rsvhj1}) to a one-dimensional invariant
manifold. This manifold is tangent at the origin to the kernel 
$\RR G$ of $\mathcal{L}$, and solutions to (\ref{rsvhj1}) which lie 
on this manifold satisfy $v(\tau,\xi) \approx M(\tau)\,G(\xi)$ for 
large times. Inserting this ansatz into (\ref{rsvhj1}) and 
integrating over $\RR^N$ we obtain the ordinary differential 
equation $dM/d\tau + \|\nabla G\|_{L^{q_\star}}^{q_\star}
\,M^{q_\star}=0$ for $M(\tau)$, from which we deduce that $M(\tau)
\approx M_\star\,\tau^{-(N+1)}$ for large times. Returning
to the original variables $(t,x)$, we then conclude that $u(t)
\approx M_\star (\ln{t})^{-(N+1)} g(t)$ as $t \to \infty$, and 
Theorem~\ref{th1} follows. 

To construct the center manifold, it is necessary to assume that 
the solutions we consider decay a little bit faster as $|x|\to
\infty$ than what is needed to be integrable. Indeed, using
the results of \cite[Appendix~A]{GW02}, it is easy to see that
the spectrum of the linear operator $\mathcal{L}$ in $L^1(\RR^N)$
is just the left half-plane $\{z \in \CC\,|\, \Re(z) \le 0\}$
(no spectral gap). In contrast, the spectrum of the same operator 
in $L^2_m(\RR^N) = L^2(\RR^N ; \left( 1+|\xi|^{2m} 
\right) \dd\xi)$ is given by 
$$
\sigma(\mathcal{L};L_m^2(\RR^N)) \,=\, \Bigl\{z \in\CC \,\Big|\,
\Re(z)\le \frac{N}{4} - \frac{m}{2} \Bigr\} \,\bigcup\,
\Bigl\{ -\frac{k}{2} \,\Big|\, k\in\NN \Bigr\}\,, 
$$
see \cite[Theorem~A.1]{GW02}. Thus, if $m > N/2$, the operator
$\mathcal{L}$ has a simple isolated eigenvalue at the origin and the
rest of the spectrum is strictly contained in the interior of the left
half-plane, a spectral configuration which allows to construct the
center manifold. This explains the choice of the weighted Lebesgue
space $L^2_m(\RR^N)$ in Theorem~\ref{th1}.  In fact, since the
nonlinearity in (\ref{vhj1}) involves the gradient of the solution
$u$, we shall rather use the corresponding Sobolev space
$H^1_m(\RR^N)$ (defined in (\ref{H1mdef}) below) in the proof.

The rest of this paper is organized as follows. In the next section,
we recall existence and uniqueness results for (\ref{rsvhj1}),
(\ref{rsvhj2}) and establish the convergence to zero of the solution
$v(\tau)$ in $H^1_m(\RR^N)$ as $\tau\to\infty$. In Section~\ref{aim}, we
study a suitable truncated version of (\ref{rsvhj1}), (\ref{rsvhj2})
to which we can apply an abstract result of \cite{CHT97} to construct
the invariant manifold. The proof of Theorem~\ref{th1} is then
performed in the final section.

\section{Global existence and convergence to zero}\label{ctz} 

In this section we briefly discuss the Cauchy problem for the 
rescaled equation (\ref{rsvhj1}) and we show that the solutions 
converge to zero as $\tau \to \infty$. We first consider initial 
data in the Lebesgue space $L^1(\RR^N) \cap L^2(\RR^N)$.

\begin{prop}\label{prb1}
Let $u_0$ be a non-negative function in $L^1(\RR^N) \cap L^2(\RR^N)$. 
Then (\ref{rsvhj1}), (\ref{rsvhj2}) have a unique non-negative 
(mild) solution  
$$
v\in \mathcal{C}([0,\infty);L^1(\RR^N) \cap L^2(\RR^N))\,\cap\,
L^\infty_{\mathrm{loc}}((0,\infty);W^{1,\infty}(\RR^N))\,, 
$$
which moreover satisfies 
\beqn
\label{b0}
\lim_{\tau\to\infty} \left\{ \|v(\tau)\|_{L^1} + \|
v(\tau)\|_{L^\infty} + \|\nabla v(\tau) \|_{L^\infty}
\right\} \,=\, 0\,.  
\eeqn 
\end{prop} 

\noindent\textbf{Proof:} For such initial data $u_0$, the results 
of \cite{BL99,BSW02,GGK03} imply that the original system 
(\ref{vhj1}), (\ref{vhj2}) has a unique (mild) solution
$$
u\in \mathcal{C}([0,\infty);L^1(\RR^N) \cap L^2(\RR^N))\,\cap\,
\mathcal{C}((0,\infty);W^{1,\infty}(\RR^N))\,. 
$$
For all $t > 0$, the function $u(t,x)$ is $\mathcal{C}^1$
with respect to $t$, $\mathcal{C}^2$ with respect to $x$, and 
(\ref{vhj1}) is satisfied in the classical sense. In addition, the 
following bounds hold for all $t > 0$:
\beqn
\label{b2}
\|u(t)\|_{L^1} + t^{N/2}\,\|u(t)\|_{L^\infty} +
t^{(N+1)/2}\,\|\nabla u(t)\|_{L^\infty} \,\le\, C\, \|u(t/2)\|_{L^1}
\,\le\, C\, \|u_0\|_{L^1}\,. 
\eeqn
Since $\|u(t)\|_{L^1}\to 0$ as $t\to\infty$ by
\cite{BL99,BK99}, we deduce from (\ref{b2}) that  
\beqn
\label{bb0}
\lim_{t\to\infty} \left\{ \|u(t)\|_{L^1} + t^{N/2}\,\|
u(t)\|_{L^\infty} + t^{(N+1)/2}\,\|\nabla u(t)\|_{L^\infty}
\right\} \,=\, 0\,. 
\eeqn
The conclusions of Proposition~\ref{prb1} are straightforward 
consequences of these results, since (\ref{rsvhj1}) is 
obtained from (\ref{vhj1}) via the simple transformation
(\ref{newuf}). In particular, $\|v(\tau)\|_{L^1} = \|u(t)\|_{L^1}$, 
$\|v(\tau)\|_{L^\infty} = (1{+}t)^{N/2}\|u(t)\|_{L^\infty}$,
$\|\nabla v(\tau)\|_{L^\infty} = (1{+}t)^{(N+1)/2}\|\nabla u(t)
\|_{L^\infty}$, hence (\ref{b0}) follows immediately from 
(\ref{bb0}). \qed 

\medskip
We next study the properties of solutions to (\ref{rsvhj1}),
(\ref{rsvhj2}) in the weighted Lebesgue space $L^2_m(\RR^N)$ 
defined in (\ref{L2mdef}), and in the corresponding Sobolev 
space 
\beqn
\label{H1mdef}
H^1_m(\RR^N) \,=\, \Bigl\{v \in H^1(\RR^N) \,\Big|\, 
\|v\|_m := (|v|_m^2 + |\nabla v|_m^2)^{1/2} < \infty\Bigr\}\,,
\eeqn
where
$$
|v |_m  \,=\, \left( \int_{\RR^N} \left( 1+|\xi|^{2m}
\right)\,|v(\xi)|^2 \dd\xi\right)^{1/2}\,.
$$

\begin{prop}\label{prb2}
Let $u_0$ be a non-negative function in $L^2_m(\RR^N)$ for some
$m > N/2$. Then the solution $v$ to (\ref{rsvhj1}), (\ref{rsvhj2}) 
given by Proposition~\ref{prb1} satisfies
$$
v \in \mathcal{C}([0,\infty);L^2_m(\RR^N)) \,\cap\,
\mathcal{C}((0,\infty);H^1_m(\RR^N))\,, 
$$
and 
\beqn
\label{b3}
\lim_{\tau\to\infty} \|v(\tau)\|_{m} \,=\, 0\,. 
\eeqn 
\end{prop} 

\noindent\textbf{Proof:} The fact that $v \in \mathcal{C}([0,T);
L^2_m(\RR^N)) \cap \mathcal{C}((0,T);H^1_m(\RR^N))$ for some 
$T > 0$ can be established by a classical fixed point argument, 
which will be implemented in Section~\ref{aim} for a truncated 
version of (\ref{rsvhj1}). Here we just obtain differential
inequalities for the norms $|v(\tau)|_m$ and $|\nabla v(\tau)|_m$ 
which imply (in view of the local existence theory) that $T > 0$ 
can be chosen arbitrarily large and that (\ref{b3}) holds. 

We first multiply (\ref{rsvhj1}) by $v$ and integrate over 
$\RR^N$. Using the non-negativity of $v$ and integrating by 
parts, we find
\beqn
\label{bb1}
\frac{1}{2}\,\frac{d}{d\tau} \int_{\RR^N} v^2\dd\xi 
\,=\, \int_{\RR^N} v \,\partial_\tau v\dd\xi 
\,\le\,  \int_{\RR^N} v \mathcal{L}v\dd\xi \,=\, 
-\int_{\RR^N} |\nabla v|^2\dd\xi + \frac{N}{4} \int_{\RR^N}
v^2\dd\xi\,.
\eeqn
Similarly, multiplying (\ref{rsvhj1}) by $|\xi|^{2m}\,v$,
we obtain
\bear
\label{bb2}
&&\frac{1}{2}\,\frac{d}{d\tau} \int_{\RR^N} |\xi|^{2m}\,v^2\
\dd\xi \,\le\,  \int_{\RR^N} |\xi|^{2m}\,v \Bigl(\Delta v + 
\frac{1}{2}\,\xi\cdot\nabla v + \frac{N}{2}v\Bigr)\dd\xi \\ \nonumber
&&\qquad \,=\, - \int_{\RR^N} |\xi|^{2m}\,|\nabla
v|^2\dd\xi  \,-\, 2m \int_{\RR^N} |\xi|^{2m-2}\,v\,
\xi \cdot \nabla v \dd\xi \,-\, \frac{2m{-}N}{4} \int_{\RR^N} 
|\xi|^{2m}\,v^2\dd\xi\,.
\eear
The only difficulty is to bound the integral involving $\xi \cdot
\nabla v$. If $m \ge 1$, we have by Young's inequality
$$
2m|\xi|^{2m-1}v|\nabla v| \,\le\, \frac12 |\xi|^{2m}|\nabla v|^2 + 
2m^2 |\xi|^{2m-2}v^2\,, \quad \hbox{and} \quad 
|\xi|^{2m-2} \,\le\, \varepsilon |\xi|^{2m} + C(\varepsilon)\,,
$$
where $\varepsilon > 0$ is arbitrary. If $1/2 < m < 1$ 
(which is possible only if $N = 1$) we find similarly
$$
2m|\xi|^{2m-1}v|\nabla v| \,\le\, \frac12 |\nabla v|^2 + 
2m^2 |\xi|^{4m-2}v^2\,, \quad \hbox{and} \quad 
|\xi|^{4m-2} \,\le\, \varepsilon |\xi|^{2m} + C(\varepsilon)\,.
$$
In both cases, summing up (\ref{bb1}) and (\ref{bb2}), we obtain 
the inequality
\beqn
\label{bb3}
\frac{d}{d\tau}\, |v(\tau)|_m^2 + |\nabla v(\tau)|_m^2 
+ \frac{2m{-}N{-}8m^2\varepsilon}{2}\, |v(\tau)|_m^2 \,\le\, 
(m+4m^2C(\varepsilon))\,\|v(\tau)\|_{L^2}^2\,.
\eeqn
We now choose $\varepsilon > 0$ sufficiently small so that 
$2m - N -8 m^2\varepsilon > 0$. Since $\|v(\tau)\|_{L^2}\to 0$ by 
(\ref{b0}), the differential inequality (\ref{bb3}) implies that 
$|v(\tau)|_{m} \to 0$ as $\tau\to \infty$. 

\medskip
We next control the evolution of $\nabla v$. Multiplying (\ref{rsvhj1})
by $-\Delta v$ and integrating by parts, we obtain
\beqn
\label{bb4}
\frac{1}{2}\,\frac{d}{d\tau} \int_{\RR^N} |\nabla v|^2 \dd\xi 
\,=\, -\int_{\RR^N} |\Delta v|^2 \dd\xi + \frac{N{+}2}{4} 
\int_{\RR^N} |\nabla v|^2 \dd\xi + \int_{\RR^N} \Delta v\, 
|\nabla v|^{q_\star}\dd\xi\,.
\eeqn
Similarly,
\bear
\nonumber
&&\frac{1}{2}\,\frac{d}{d\tau} \int_{\RR^N} |\xi|^{2m}\,|\nabla v|^2 
\dd\xi \,=\, -\int_{\RR^N} |\xi|^{2m}\,|\Delta v|^2 \dd\xi + 
\frac{N{+}2-2m}{4} \int_{\RR^N} |\xi|^{2m}\,|\nabla v|^2 \dd\xi\\
\label{bb5}
&& -2m \int_{\RR^N} |\xi|^{2m-2}\,\Delta v\,\xi\cdot\nabla v
\dd\xi + \int_{\RR^N} \Bigl(|\xi|^{2m}\,\Delta v + 2m |\xi|^{2m-2}\,
\xi\cdot\nabla v\Bigr)|\nabla v|^{q_\star}\dd\xi\,.
\eear
Using the crude estimate $|\xi|^{2m-1} \le 1+|\xi|^{2m}$, we find
$$
2m \int_{\RR^N} |\xi|^{2m-1}\,|\Delta v|\,|\nabla v| \dd\xi
\,\le\, \frac14\,|\Delta v|_m^2 + C |\nabla v|_m^2\,.
$$
Moreover, as $\|\nabla v\|_{L^\infty}$ is uniformly bounded for all
$\tau \ge 1$ by (\ref{b0}), we have for such times
\bean
\int_{\RR^N} (1{+}|\xi|^{2m})\,|\Delta v|\,|\nabla v|^{q_\star}
\dd\xi &\le& \frac14\,|\Delta v|_m^2 + C |\nabla v|_m^2\,, \\
2m \int_{\RR^N} |\xi|^{2m-1} \,|\nabla v|\,|\nabla v|^{q_\star}
\dd\xi &\le& C |\nabla v|_m^2\,.
\eean
Thus adding up (\ref{bb4}) and (\ref{bb5}), we obtain
\beqn
\label{bb6}
\frac{d}{d\tau}\, |\nabla v(\tau)|_m^2 + |\Delta v(\tau)|_m^2 
+ |\nabla v(\tau)|_m^2 \,\le\, K\,|\nabla v(\tau)|_m^2\,, 
\quad \tau \ge 1\,,
\eeqn
for some $K > 0$. Now, if we combine (\ref{bb3}) and (\ref{bb6}), 
we see that $h(\tau) := K|v(\tau)|_m^2 + |\nabla v(\tau)|_m^2$ 
satisfies a differential inequality of the form
$$
  h'(\tau) + \varepsilon_0\,h(\tau) \,\le\, C\,\|v(\tau)\|_{L^2}^2\,,
  \quad \tau \ge 1\,,
$$
for some positive constants $\varepsilon_0$ and $C$. Thus 
$h(\tau) \to 0$ as $\tau \to \infty$, and (\ref{b3}) follows. 
This concludes the proof of Proposition~\ref{prb2}. \qed

\section{Construction of the center manifold}\label{aim}

In this section, we proceed along the lines of \cite[Section~3]{GW02}
to describe the large time behavior of the non-negative solutions to
(\ref{rsvhj1}), (\ref{rsvhj2}) in the space $L^2_m(\RR^N)$.  By
Proposition~\ref{prb2} these solutions converge to zero in
$H^1_m(\RR^N)$ as $\tau\to\infty$, hence the large time asymptotics
remain unchanged if we truncate the nonlinearity in (\ref{rsvhj1})
outside a neighborhood of the origin. This modification will 
allow us to apply the center manifold theorem as stated in 
\cite{CHT97}. 

Throughout this section, we fix a function $\chi\in\mathcal{C}^\infty
([0,\infty))$ such that $0\le\chi\le 1$, $\chi(r) = 0$ if $r\ge 4$ and
$\chi(r) = 1$ if $r\le 1$. For $\varrho>0$ and $r\ge 0$, we denote 
$\chi_\varrho(r) =\chi(r/\varrho^2)$. Given $\varrho\in (0,1)$
and $m > N/2$, we consider the initial-value problem
\bear
\label{c2}
\partial_\tau v & = & \mathcal{L} v - F_\varrho(v)\,, \qquad
(\tau,\xi)\in (0,\infty)\times\RR^N\,, \\ 
\label{c3}
v(0) & = & v_0\in H^1_m(\RR^N)\,, \quad \xi\in\RR^N\,,
\eear
where $\mathcal{L}$ is the linear operator (\ref{opL}) and
$F_\varrho$ is the truncated nonlinearity
\beqn
\label{c4}
F_\varrho(v) \,=\, \chi_\varrho\left( \|v\|_m^2 \right)\
|\nabla v|^{q_\star}\,, \quad v\in H^1_m(\RR^N)\,. 
\eeqn 
We first establish the well-posedness of (\ref{c2}), (\ref{c3}) and
show that this system generates a $\mathcal{C}^1$-smooth semiflow in
$H^1_m(\RR^N)$.  

\begin{prop}\label{prc1}
Fix $\varrho \in (0,1)$ and $m > N/2$. For each $v_0\in H^1_m(\RR^N)$,
the initial-value problem (\ref{c2}), (\ref{c3}) has a unique global
solution $v\in \mathcal{C}([0,\infty);H^1_m(\RR^N))$. Moreover, 
the map $\varphi_\tau : H^1_m(\RR^N) \to H^1_m(\RR^N)$ defined 
for $\tau \ge 0$ by $\varphi_\tau(v_0)=v(\tau)$ is globally Lipschitz
continuous, uniformly in $\tau$ on compact intervals. Finally
$\varphi_\tau$ is $\mathcal{C}^1$-smooth for all $\tau \ge 0$, so that
the family $(\varphi_\tau)_{\tau\ge 0}$ defines a $\mathcal{C}^1$
semiflow in $H^1_m(\RR^N)$.
\end{prop}

Before proving Proposition~\ref{prc1}, we recall that the linear
operator $\mathcal{L}$ defined by (\ref{opL}) is the generator of a
strongly continuous semigroup $\left( e^{\tau\mathcal{L}}
\right)_{\tau\ge 0}$ in $L^2_m(\RR^N)$, see e.g. \cite[Appendix~A]{GW02}.
If $m > N/2$, this semigroup is uniformly bounded, i.e. there exists
$C_1 > 0$ such that, for all $w \in L^2_m(\RR^N)$,
\beqn
\label{c5}
\left|e^{\tau\mathcal{L}} w \right|_m \,\le\, C_1\,|w|_m\,, \quad
\tau \ge 0\,.
\eeqn
More generally, let $b(\xi) = (1+|\xi|^2)^{1/2}$ and assume that 
$b^m w \in L^p(\RR^N)$ for some $p \in [1,2]$. Then $e^{\tau\mathcal{L}}
w \in L^2_m(\RR^N)$ for all $\tau > 0$, and there exists $C_2 > 0$
such that
\beqn
\label{c51}
\left|e^{\tau\mathcal{L}} w \right|_m \,\le\, \frac{C_2}
{a(\tau)^{\frac{N}{2}(\frac1p-\frac12)}}\,\|b^m w\|_{L^p}\,,
\quad \tau > 0\,,
\eeqn
where $a(\tau) = 1 - e^{-\tau}$, see \cite[Proposition~A.5]{GW02} and
\cite[Proposition~A.2]{GW02}. Similar estimates hold for the spatial 
derivatives of $e^{\tau\mathcal{L}} w$. For instance, as $\nabla
e^{\tau\mathcal{L}} = e^{\tau/2}\,e^{\tau\mathcal{L}} \nabla$, 
it follows from (\ref{c5}) that $|\nabla e^{\tau\mathcal{L}} w|_m 
\le C_1 \,e^{\tau/2}|\nabla w|_m$ for all $w \in H^1_m(\RR^N)$. 
In addition, if $b^m w \in L^p(\RR^N)$ for some $p \in [1,2]$, 
we have the analog of (\ref{c51}):
\beqn
\label{c52}
\left|\nabla e^{\tau\mathcal{L}} w \right|_m \,\le\, \frac{C_3}
{a(\tau)^{\frac{N}{2}(\frac1p-\frac12)+\frac12}}\,\|b^m w\|_{L^p}\,,
\quad \tau > 0\,.
\eeqn

In the rest of this section, we fix $p\in (1,2)$ such that
\beqn
\label{c7}
  \frac{2(N+1)}{N+3} \,<\, p \,<\, \frac{2(N+1)}{N+2}\,.
\eeqn
Given $T>0$ and $w\in\mathcal{C}([0,T];H_m^1(\RR^N))$, we define 
\beqn
\label{c6}
(\mathcal{N}_\varrho w)(\tau) \,=\, \int_0^\tau e^{(\tau-s)\mathcal{L}}\,
F_\varrho(w(s))\dd s\,, \quad \tau\in [0,T]\,. 
\eeqn
Then $\mathcal{N}_\varrho w$ belongs to $\mathcal{C}([0,T];H_m^1(\RR^N))$ 
and enjoys the following property: 

\begin{lem}\label{lec2} There exists a constant $C_4 > 0$ such that,
for all $T > 0$, all $\varrho \in (0,1)$, and all $w_1, w_2 \in \mathcal{C}
([0,T];H_m^1(\RR^N))$, the following inequality holds:
$$
  \left\|\left( \mathcal{N}_\varrho w_1 - \mathcal{N}_\varrho w_2 
  \right)(\tau)\right\|_m \,\le\, C_4\,Z_p(\tau)\,\varrho^{q_\star-1}\,
  \sup_{s\in [0,\tau]}\|(w_1-w_2)(s)\|_m\,, \quad \tau\in [0,T]\,,
$$
where
$$
Z_p(\tau) \,=\, \int_0^\tau \frac{1}{a(s)^{\frac{N}{2}(\frac1p-\frac12)}}
\Bigl(1 + \frac{1}{a(s)^{1/2}}\Bigr)\dd s\,.
$$
\end{lem}

\medskip\noindent\textbf{Proof:} We first observe that our choice of 
$p$ in (\ref{c7}) guarantees that $\frac{N}{2}(\frac1p-\frac12) < 
\frac12$, so that $Z_p(\tau)$ is well-defined and finite for every 
$\tau\ge 0$. The following inequality will also be useful: if 
$f,g \in L^2_m(\RR^N)$, then $b^m |f|^{q_\star-1} g \in L^p(\RR^N)$ 
and 
\beqn
\label{lec3}
\left\|b^m\,|f|^{q_\star-1}\,g\right\|_{L^p} \,\le\, C_5\
|f|_m^{q_\star-1}\,|g|_m\,. 
\eeqn
Indeed, by H\"older's inequality, 
$$
\|b^m |f|^{q_\star-1} g\|_{L^p} \,\le\, \|b^m g\|_{L^2} 
\|f\|_{L^r}^{q_\star-1}\,, \quad \hbox{where}\quad 
r \,=\, (q_\star-1)\,\frac{2p}{2-p}\,.
$$
But $\|b^m g\|_{L^2} \le C |g|_m$, and $1 < r < 2$ by (\ref{c7}), 
hence $\|f\|_{L^r} \le C (\|f\|_{L^1}+\|f\|_{L^2}) \le C|f|_m$
because $L^2_m(\RR^N) \hookrightarrow L^1(\RR^N) \cap L^2(\RR^N)$
for $m > N/2$. This proves (\ref{lec3}). 

Now, let $T > 0$, $\varrho \in (0,1)$, and $w_1, w_2 \in \mathcal{C}
([0,T];H_m^1(\RR^N))$. For all $\tau \in [0,T]$ we have by 
(\ref{c6}), (\ref{c51})
\bear
\nonumber
\left|\left( \mathcal{N}_\varrho w_1 - \mathcal{N}_\varrho w_2 \right)
(\tau)\right|_m & \le & \int_0^\tau \left|e^{(\tau-s)\mathcal{L}} 
\left(F_\varrho(w_1(s)) - F_\varrho(w_2(s))\right) \right|_m
\dd s \\ \label{c8}
& \le & \int_0^\tau \frac{C_2}{a(\tau{-}s)^{\frac{N}{2}(\frac1p-\frac12)}} 
\,\left\|b^m \left(F_\varrho(w_1(s)) - F_\varrho(w_2(s)) \right) 
\right\|_{L^p}\dd s\,.
\eear
For a fixed $s \in [0,\tau]$, we can assume for instance that 
$\|w_1(s)\|_m \ge \|w_2(s)\|_m$. Then 
\bean
\left\|b^m \left( F_\varrho(w_1(s)) - F_\varrho(w_2(s)) \right)
\right\|_{L^p} & \le & \left|\chi_\varrho\left( \|
w_1(s)\|_m^2 \right) - \chi_\varrho\left( \|w_2(s)\|_m^2
\right) \right|\, \left\|b^m\,|\nabla
w_2(s)|^{q_\star} \right\|_{L^p} \\ 
& + & \chi_\varrho\left( \|w_1(s)\|_m^2 \right)\,\left\|
b^m\,\left( |\nabla w_1(s)|^{q_\star} - |\nabla
w_2(s)|^{q_\star} \right) \right\|_{L^p}\,. 
\eean
Obviously, the right-hand side vanishes if $\|w_2(s)\|_m \ge 2\varrho$, 
hence we can suppose that $\|w_2(s)\|_m \le 2\varrho$. To bound
the first term, we apply (\ref{lec3}) with $f = g = |\nabla w_2|$ and 
obtain $\|b^m\,|\nabla w_2(s)|^{q_\star}\|_{L^p} \le C_5 |\nabla 
w_2(s)|_m^{q_\star} \le C \varrho^{q_\star}$. Moreover, if 
$\|w_1(s)\|_m \le 4\varrho$, we have
$$
\left|\chi_\varrho\left(\|w_1(s)\|_m^2 \right) -
\chi_\varrho\left(\|w_2(s)\|_m^2 \right) \right|
\,\le\, \frac{C}{\varrho^2} \left(\|w_1(s)\|_m^2 - \|w_2(s)\|_m^2
\right) \,\le\, \frac{C}{\varrho}\|(w_1-w_2)(s)\|_m\,,
$$
and the same estimate holds if $\|w_1(s)\|_m \ge 4\varrho$ because
$\|(w_1-w_2)(s)\|_m \ge 2\varrho$ in that case. Thus
$$
\left|\chi_\varrho\left(\|w_1(s)\|_m^2 \right) -
\chi_\varrho\left(\|w_2(s)\|_m^2 \right) \right|\
\left\|b^m\,|\nabla w_2(s)|^{q_\star} \right\|_{L^p}
\,\le\, C\,\varrho^{q_\star-1}\,\|(w_1-w_2)(s)\|_m\,.
$$
On the other hand, using (\ref{lec3}) and the inequality 
$|\,|y|^{q_\star}-|z|^{q_\star}\,| \le {q_\star}(|y|^{q_\star-1}+
|z|^{q_\star-1})|y-z|$, we obtain
\bean
& & \chi_\varrho\left( \|w_1(s)\|_m^2 \right)\,\left\|b^m
\left( |\nabla w_1(s)|^{q_\star} - |\nabla
w_2(s)|^{q_\star} \right) \right\|_{L^p} \\ 
&& \qquad \le\, C\,\chi_\varrho\left( \|w_1(s)\|_m^2 \right)\,\left(
\left|\nabla w_1(s) \right|_m^{q_\star-1} + \left|\nabla w_2(s)
\right|_m^{q_\star-1} \right)\,\left|\nabla(w_1-w_2)(s)
\right|_m \\ 
&& \qquad \le\, C\,\chi_\varrho\left( \|w_1(s)\|_m^2 \right)\,\|
w_1(s)\|_m^{q_\star-1}\,\left\|(w_1-w_2)(s) \right\|_m \\ 
&& \qquad \le\, C\,\varrho^{q_\star-1}\,\|(w_1-w_2)(s)\|_m\,. 
\eean
Therefore $\|b^m (F_\varrho(w_1(s)) - F_\varrho(w_2(s)))\|_{L^p} \le 
C\,\varrho^{q_\star-1} \|(w_1-w_2)(s)\|_m$, and inserting this bound
into (\ref{c8}) we conclude that
\beqn
\label{c9}
\left|\left( \mathcal{N}_\varrho w_1 - \mathcal{N}_\varrho w_2 
\right)(\tau) \right|_m \,\le\, C\,\varrho^{q_\star-1}\,\int_0^\tau
\frac{1}{a(\tau{-}s)^{\frac{N}{2}(\frac1p-\frac12)}}\,
\left\|(w_1-w_2)(s) \right\|_m\dd s\,.
\eeqn
Finally, using (\ref{c52}) and proceeding in the same way, we also 
obtain
\bear
\nonumber
\left|\nabla \left(\mathcal{N}_\varrho w_1 - \mathcal{N}_\varrho w_2 \right)
(\tau)\right|_m & \le & \int_0^\tau \frac{C_3}{a(\tau{-}s)^{\frac{N}{2}
(\frac1p-\frac12)+\frac12}} \,\left\|b^m \left(F_\varrho(w_1(s)) - 
F_\varrho(w_2(s)) \right) \right\|_{L^p}\dd s \\ \label{c10}
& \le & C\,\varrho^{q_\star-1} \int_0^\tau
\frac{1}{a(\tau{-}s)^{\frac{N}{2}(\frac1p-\frac12)+\frac12}}\,
\left\|(w_1-w_2)(s) \right\|_m\dd s\,.
\eear
Lemma~\ref{lec2} is now a immediate consequence of (\ref{c9}) and 
(\ref{c10}). \qed 

\medskip
\noindent\textbf{Proof of Proposition~\ref{prc1}:} Given $v_0\in
H^1_m(\RR^N)$, we choose $K > 2C_1\|v_0\|_m$ and $T > 0$ sufficiently 
small so that
\beqn
\label{contrac}
  C_1 \,\|v_0\|_m\,e^{T/2} \,\le\, \frac{K}{2}\,, \quad 
  \hbox{and}\quad C_4\, \varrho^{q_\star-1} Z_p(T) \,\le\, 
  \frac12\,, 
\eeqn
where $C_1$ is as in (\ref{c5}) and $C_4, Z_p$ are defined 
in Lemma~\ref{lec2}. We introduce the set 
$$
X_{K,T} \,=\, \Bigl\{ w\in\mathcal{C}([0,T];H^1_m(\RR^N))~\Big|~
\sup_{\tau\in [0,T]} \|w(\tau)\|_m\le K \Bigr\}\,, 
$$
which is a complete metric space for the distance $d_T$ defined by
$$
d_T(w_1,w_2) \,=\, \sup_{\tau\in [0,T]} \|
(w_1-w_2)(\tau)\|_m\,,\quad (w_1,w_2)\in X_{K,T}\times X_{K,T}\,. 
$$
Using (\ref{c5}) and Lemma~\ref{lec2} it is straightforward 
to verify that, if $w\in X_{K,T}$, then the function 
$\mathcal{T}_\varrho w : [0,T] \to H^1_m(\RR^N)$ defined by  
$$
(\mathcal{T}_\varrho w)(\tau) \,=\, e^{\tau\mathcal{L}} v_0 - 
(\mathcal{N}_\varrho w)(\tau)\,,\quad \tau\in [0,T]\,, 
$$ 
belongs to $X_{K,T}$, and that the map $w \mapsto \mathcal{T}_\varrho w$ 
is a strict contraction in $X_{K,T}$. By the Banach fixed point theorem,
$\mathcal{T}_\varrho$ has thus a unique fixed point $v$ in $X_{K,T}$. 
This proves that the Cauchy problem (\ref{c2}), (\ref{c3}) is locally
well-posed in $H^1_m(\RR^N)$. 

Let $T_*(v_0) \in (0,\infty]$ be the maximal existence time
for the solution of (\ref{c2}), (\ref{c3}) in $H^1_m(\RR^N)$.   
For all $\tau < T_*(v_0)$, it follows from (\ref{c5}), (\ref{c9}), 
and (\ref{c10}) (with $w_1=v$ and $w_2=0$) that 
$$
\|v(\tau)\|_m \,\le\, C_1\,e^{\tau/2}\,\|v_0\|_m + C_4\
\varrho^{q_\star-1} \int_0^\tau  \frac{\left\|v(s)
\right\|_m}{a(\tau{-}s)^{\frac{N}{2}(\frac1p-\frac12)+\frac12}} \dd s\,. 
$$
Using a version of Gronwall's lemma (see e.g. \cite[Lemma~7.1.1]{He81}), 
we deduce that $\|v(\tau)\|_m$ cannot blow  up in finite time, 
hence $T_*(v_0) = \infty$. Thus (\ref{c2}) has a unique global 
solution $v\in \mathcal{C}([0,\infty);H^1_m(\RR^N))$ for all 
$v_0 \in H^1_m(\RR^N)$, and we may define a semiflow 
$(\varphi_\tau)_{\tau\ge 0}$ by the relation
$\varphi_\tau(v_0)=v(\tau)$ for $\tau\ge 0$. 

By construction, the map $v_0 \mapsto \varphi_\tau(v_0)$ is 
globally Lipschitz continuous, uniformly in time on compact intervals:
for each $T > 0$, there exists $L(T) > 0$ such that
\beqn
\label{corc5}
\left\|\varphi_\tau(v_0) - \varphi_\tau\left( \hat{v}_0 \right)
\right\|_m \,\le\, L(T)\,\left\|v_0 - \hat{v}_0 \right\|_m\,, 
\eeqn
for all $\tau\in [0,T]$ and all $(v_0,\hat{v}_0)\in H^1_m(\RR^N)\times
H^1_m(\RR^N)$. Indeed, by the semigroup property, it is sufficient to
prove (\ref{corc5}) for a $T > 0$ satisfying (\ref{contrac}), in which
case (\ref{corc5}) follows immediately from the fixed point argument
above, with $L(T) = 2C_1\,e^{T/2}$. This proof also shows that $L(T)$
can be chosen independent of $\varrho$ if $\varrho \in (0,1)$.
Finally, the fact that the map $\varphi_\tau$ is $\mathcal{C}^1$ for
each $\tau\ge 0$ is obtained by classical arguments which we omit
here.  We only mention that, given $v_0\in H^1_m(\RR^N)$, $\tau\ge 0$,
and $h\in H^1_m(\RR^N)$, the differential $D\varphi_\tau(v_0)h$ of
$\varphi_\tau$ at $v_0$ applied to $h$ is equal to $V(\tau)$, where
$V$ denotes the solution of the linear non-autonomous equation
\bean
\partial_\tau V &=& \mathcal{L} V - q_\star \chi_\varrho\left( \|
v\|_m^2 \right)\,|\nabla v|^{q_\star-2}\,\nabla v \cdot
\nabla V - 2 \chi_\varrho'\left( \|v\|_m^2 \right)\
|\nabla v|^{q_\star}\,\ll v , V \gg_m\,, \\
V(0) &=& h\,. 
\eean
Here $v(\tau) = \varphi_\tau(v_0)$ and $\ll \cdot,\cdot \gg_m$ denotes 
the scalar product in $H^1_m(\RR^N)$. In particular, since 
$\varphi_\tau(0)=0$ for all $\tau\ge 0$, this formula shows that
$D\varphi_\tau(0)=e^{\tau\mathcal{L}}$ for each $\tau\ge 0$. \qed 

\begin{rem}\label{rec4}
It can actually be shown that the differential $D\varphi_\tau : 
H^1_m(\RR^N) \to \mathscr{L}(H^1_m(\RR^N))$ is H\"older continuous 
with exponent $q_\star-1$ for any $\tau \ge 0$. 
\end{rem}

For later use, we also point out the following properties of the 
time-one map $\varphi_1$:

\begin{cor}\label{lec6}
The map $\mathcal{R} = \varphi_1 - e^\mathcal{L}$ belongs to 
$\mathcal{C}^1(H^1_m(\RR^N);H^1_m(\RR^N))$ and satisfies 
$\mathcal{R}(0)=0$, $D\mathcal{R}(0)=0$. Moreover $\mathcal{R}$
is globally Lipschitz continuous and there exists $C_6 > 0$
(independent of $\varrho$) such that its Lipschitz constant satisfies 
$\mathrm{Lip}(\mathcal{R}) \le C_6\,\varrho^{q_\star-1}$. 
\end{cor}

\noindent\textbf{Proof:} We know from Proposition~\ref{prc1} that 
$\mathcal{R}$ is indeed a $\mathcal{C}^1$-map from $H^1_m(\RR^N)$ 
into itself, and it was observed at the end of the proof that
$\varphi_1(0)=0$ and $D\varphi_1(0) = e^{\mathcal{L}}$, hence
$\mathcal{R}(0)=0$ and $D\mathcal{R}(0)=0$. Now, given $v_0,\hat{v}_0$ 
in $H^1_m(\RR^N)$ we define $v(\tau) =\varphi_\tau(v_0)$ and
$\hat{v}(\tau) = \varphi_\tau\left(\hat{v}_0 \right)$ for $\tau\ge
0$. Using Lemma~\ref{lec2} and estimate (\ref{corc5}) we find
\bean
\left\|\mathcal{R}(v_0) - \mathcal{R}\left( \hat{v}_0 \right)
\right\|_m  & = & \left\|(\mathcal{N}_\varrho v)(1) -
(\mathcal{N}_\varrho \hat{v})(1) \right\|_m \,\le\, 
C_4\,\varrho^{q_\star-1}\,\sup_{s\in [0,1]} \left\|\left( v - \hat{v}
\right)(s)\right\|_m \\ 
& \le & C_4\,L(1)\,\varrho^{q_\star-1}\,\left\|v_0 - \hat{v}_0
\right\|_m\,, 
\eean
which is the desired bound. \qed

\medskip
Having associated a $\mathcal{C}^1$-semiflow to the truncated
system (\ref{c2}),  we now turn to the construction of a 
center manifold for this semiflow at the origin. If $m > N/2$, 
we can decompose $H^1_m(\RR^N) = E_c \oplus E_s$, where $E_c = 
\{\alpha G\,|\,\alpha \in \RR\}$ is the kernel of the operator 
$\mathcal{L}$ and 
\beqn
\label{c11}
E_s \,=\, \Bigl\{w\in H^1_m(\RR^N) ~\Big|~
\int_{\RR^N} w(\xi)\dd\xi = 0 \Bigr\}\,.
\eeqn 
We recall that $G$ is the Gaussian function defined in (\ref{gk}). 
Let $P_0$ be the continuous projection onto $E_c$ along $E_s$,
namely
  $$
  P_0 w \,=\, \left(\int_{\RR^N} w(\xi)\dd\xi \right)\,G\,, 
  \quad w\in H^1_m(\RR^N)\,,
$$
and let $Q_0 = \mathbf{1} - P_0$. It is easily verified that
$P_0$ and $Q_0$ commute with $\mathcal{L}$, so that the subspaces
$E_c$ and $E_s$ are invariant under the action of $\mathcal{L}$. 
Moreover, we know from \cite[Appendix~A]{GW02} that the spectrum 
of the restriction of $\mathcal{L}$ to the invariant subspace 
$E_s$ is strictly contained in the left-half plane in $\CC$, 
because the associated semigroup $e^{\tau\mathcal{L}}$ decreases 
exponentially in $E_s$. More precisely, if $\mu_0\in (0,1/2)$ 
satisfies $2\mu_0 < m-(N/2)$, there exists $C_7 > 0$ such that 
\beqn
\label{expdecay}
\left|e^{\tau\mathcal{L}} Q_0 w \right|_m + a(\tau)^{1/2}\,
\left|\nabla e^{\tau\mathcal{L}} Q_0 w \right|_m \,\le\, C_7
\,e^{-\mu_0 \tau}\,|w|_m\,, 
\eeqn
for all $w\in L^2_m(\RR^N)$ and all $\tau > 0$, see 
\cite[Proposition~A.2]{GW02}. 

We are now in a position to apply the invariant manifold theorem 
as stated in \cite[Theorem~1.1]{CHT97}. The main result of 
this section reads: 

\begin{theo}\label{thc7}
Fix $\mu\in (0,1/2)$ such that $2\mu < m-(N/2)$. If $\varrho > 0$ is 
sufficiently small, there exists a globally Lipschitz continuous map
$f\in\mathcal{C}^1(E_c;E_s)$ with $f(0)=0$ and $Df(0)=0$ such that 
the submanifold $W_c = \{\alpha G + f(\alpha G) \,|\, \alpha\in\RR\} 
\subset H^1_m(\RR^N)$ enjoys the following properties: 
\begin{enumerate}
\item[{\bf (a)}] $\varphi_\tau(W_c) = W_c$ for every $\tau\ge 0$,
\item[{\bf (b)}] for every $v_0\in H^1_m(\RR^N)$ there exist a unique
$w_0\in W_c$ and a positive constant $C_8(v_0)$ such that 
\beqn
\label{conv}
\left\|\varphi_\tau(v_0) - \varphi_\tau(w_0) \right\|_m
\,\le\, C_8(v_0)\,e^{-\mu\tau} \;\;\mbox{ for }\;\; \tau\ge 0\,. 
\eeqn
\end{enumerate}
\end{theo}

\noindent\textbf{Proof:} Theorem~\ref{thc7} readily follows from
\cite[Theorem~1.1]{CHT97} once we have checked that the assumptions
(\textbf{H.1})--(\textbf{H.4}) of \cite{CHT97} are fulfilled in our
case. By Proposition~\ref{prc1}, $(\varphi_\tau)_{\tau\ge 0}$ is a
$\mathcal{C}^1$ semiflow in $H^1_m(\RR^N)$ and $\varphi_\tau$ is
globally Lipschitz continuous, uniformly for $\tau \in [0,1]$. Therefore,
\cite[(\textbf{H.1})]{CHT97} is verified. Next, assumption
\cite[(\textbf{H.2})]{CHT97} is nothing but the decomposition
$\varphi_1=e^\mathcal{L}+\mathcal{R}$ described in
Corollary~\ref{lec6}. To check \cite[(\textbf{H.3})]{CHT97}, we 
remark that $P_0 \,e^\mathcal{L}\, P_0 = \mathbf{1}$, hence
$$
\left\|\left( P_0 \,e^\mathcal{L}\,P_0 \right)^{-k} P_0
\right\|_{\mathscr{L}(E_c)} \,=\, \|P_0\|_{\mathscr{L}(E_c)}\,,  
\quad \hbox{for all } k\in\NN\,. 
$$
On the other hand, if we choose $\mu_0 \in (\mu,1/2)$ such that
$2\mu_0 < m-(N/2)$, it follows from (\ref{expdecay}) that 
$\|e^{k\mathcal{L}} Q_0 w\|_m \le C\,e^{-k\mu_0}\,\|w\|_m$ for 
all $k \in \NN$. Since $Q_0$ and $e^\mathcal{L}$ commute, this
inequality is equivalent to
$$
\left\|\left( Q_0 \,e^\mathcal{L}\,Q_0 \right)^k Q_0
\right\|_{\mathscr{L}(E_s)} \,\le\, C\,e^{-k \mu_0}\,, \quad 
\hbox{for all } k \in\NN\,. 
$$
As $e^{-\mu_0}<1$, we have thus checked that 
\cite[(\textbf{H.3})]{CHT97} is fulfilled. Finally 
\cite[(\textbf{H.4})]{CHT97} is automatically satisfied 
if the Lipschitz constant of $\mathcal{R}$ is sufficiently small. 
By Corollary~\ref{lec6}, this is easily achieved by choosing $\varrho$ 
appropriately small.

Therefore, by \cite[Theorem~1.1]{CHT97}, there exist $\mu_1 \in 
(0,\mu_0)$ and a globally Lipschitz continuous map
$f\in\mathcal{C}^1(E_c;E_s)$ such that the submanifold  
$$
W_c \,=\, \left\{\alpha G + f(\alpha G)\,|\,\alpha\in\RR
\right\} \subset H^1_m(\RR^N)
$$
enjoys the following properties:

\medskip\noindent{\bf Invariance:} $\varphi_\tau(W_c) = W_c$ for all 
$\tau\ge 0$, and the restriction to $W_c$ of the semiflow 
$(\varphi_\tau)_{\tau\ge 0}$ can be extended to a Lipschitz continuous
flow on $W_c$.

\medskip\noindent{\bf Invariant foliation:} There is a continuous map
$h:H^1_m(\RR^N)\times E_s \to E_c$ such that, for each $v\in W_c$,
one has $h(v,Q_0 v)=P_0 v$ and the manifold 
$$
\mathcal{M}_v \,=\, \left\{ h(v,w)+w\,|\, w\in E_s \right\}
\subset H^1_m(\RR^N)
$$
passing through $v$ satisfies $\varphi_\tau(\mathcal{M}_v)\subset 
\mathcal{M}_{\varphi_\tau(v)}$ for $\tau\ge 0$ and is 
characterized by 
$$
\mathcal{M}_v \,=\, \Bigl\{ w\in H^1_m(\RR^N)~\Big|~ 
\limsup_{\tau\to\infty} \frac{1}{\tau}\, \ln{\left( \left\|
\varphi_\tau(w) - \varphi_\tau(v) \right\|_m \right)} \le - \mu_1
\Bigr\}\,. 
$$

\medskip\noindent{\bf Completeness:} For every $v\in W_c$, $\mathcal{M}_v 
\cap W_c = \{v\}$. In particular, $\mathcal{M}_v\cap \mathcal{M}_w
= \emptyset$ if $(v,w)\in W_c\times W_c$ and $v\ne w$, and  
$H^1_m(\RR^N) = \cup_{v\in W_c} \mathcal{M}_v$. 

\medskip\noindent
Moreover, we can assume that $\mu_1 \in (\mu,\mu_0)$ if $\varrho > 0$
is sufficiently small. 

\medskip
We can now conclude the proof of Theorem~\ref{thc7}. 
Assertion~\textbf{(a)} is nothing but the invariance property of 
$W_c$. To prove \textbf{(b)}, let $v_0\in H^1_m(\RR^N)$. By the
completeness property of $W_c$, there is a unique $w_0\in W_c$ 
such that $v_0\in\mathcal{M}_{w_0}$. Since $\mu < \mu_1$, we deduce 
from the invariant foliation property of $W_c$ that there is 
$\tau_0 > 0$ such that 
$$
\left\|\varphi_\tau(v_0) - \varphi_\tau(w_0) \right\|_m \,\le\, 
e^{-\mu\tau}\,, \quad \hbox{for all }\tau\ge\tau_0\,.
$$
Using in addition (\ref{corc5}), we obtain (\ref{conv}). \qed

\section{Large time behavior}\label{ltb}

This final section is entirely devoted to the proof of 
Theorem~\ref{th1}. Assume that $u_0$ is a non-negative function in 
$L^2_m(\RR^N)$, $m > N/2$, such that $\Vert u_0\Vert_{L^1} > 0$. 
Let $u(t,x)$ be the corresponding solution of (\ref{vhj1}), (\ref{vhj2}) 
and $v(\tau,\xi)$ the corresponding solution of (\ref{rsvhj1}), 
(\ref{rsvhj2}). By the strong maximum principle \cite[Corollary~
4.2]{GGK03}, we know that 
$u(t,x) > 0$ for all $t > 0$ and all $x \in \RR^N$. Choose 
$\mu \in (0,1/2)$ such that $2\mu < m-(N/2)$ and $\varrho \in (0,1)$
sufficiently small so that Theorem~\ref{thc7} applies. 

By Proposition~\ref{prb2}, the solution $v$ of (\ref{rsvhj1}) 
converges to zero in $H^1_m(\RR^N)$ as $\tau \to \infty$, hence
there exists $\tau_0 \ge 0$ such that $\|v(\tau)\|_m \le \varrho$ 
for all $\tau \ge \tau_0$. Setting $v_0 = v(\tau_0)$ and 
$\hat v(\tau) = v(\tau+\tau_0)$, we obtain a solution 
$\hat v(\tau)$ of (\ref{rsvhj1}) with initial condition $v_0 \in 
H^1_m(\RR^N)$ which satisfies 
\beqn
\label{d1}
\|\hat v(\tau)\|_m \,\le\, \varrho \quad \mbox{for all }\tau\ge 0\,. 
\eeqn
Using the notations of Section~\ref{aim}, it follows that
$\hat v(\tau) = \varphi_\tau(v_0)$ for $\tau\ge 0$, because 
(\ref{d1}) implies that $\chi_\varrho\left(\|\hat v(\tau)\|_m^2
\right) = 1$. Thus, in view of Theorem~\ref{thc7}, there exist 
$w_0 \in W_c$ and $C_9 > 0$ such that 
\beqn
\label{d2}
\left\|\hat v(\tau) - \varphi_\tau(w_0) \right\|_m \,\le\, C_9\
e^{-\mu \tau}\,, \quad \tau\ge 0\,. 
\eeqn
To simplify the notations, we set $w(\tau) = \varphi_\tau(w_0)$ 
and 
$$
  M(\tau) \,=\, \int_{\RR^N} w(\tau,\xi)\dd\xi\,, \quad \tau\ge 0\,. 
$$

We claim that 
\beqn
\label{d3}
M(\tau) > 0 \quad \mbox{for all } \tau\ge 0\,, \quad \mbox{and}\quad
\lim_{\tau\to\infty} M(\tau) = 0\,. 
\eeqn
Indeed, since $H^1_m(\RR^N) \hookrightarrow L^1(\RR^N)$, it follows 
from (\ref{d2}) that 
\beqn
\label{d4}
\left|\int_{\RR^N} \hat v(\tau,\xi)\dd\xi - M(\tau) \right| \,\le\,  
C\,\|\hat v(\tau)-w(\tau)\|_m \,\le\, C_{10}\, e^{-\mu \tau}\,, 
\eeqn
for all $\tau\ge 0$. Assume by contradiction that there exists
$\tau_1 \ge 0$ such that  $M(\tau_1) \le 0$. Since $w$ is a 
solution of (\ref{c2}), (\ref{c3}) and $F_\varrho\ge 0$, it is 
clear that $\tau\mapsto M(\tau)$ is non-increasing, hence 
$M(\tau) \le M(\tau_1)\le 0$ for $\tau \ge \tau_1$. Using (\ref{d4})
and recalling that $\hat v$ is non-negative, we thus find
$$
\|\hat v(\tau)\|_{L^1} \,=\, \int_{\RR^N} \hat v(\tau,\xi)\dd\xi 
\,\le\, M(\tau) + C_{10}\,e^{-\mu \tau} \,\le\, C_{10}\, e^{-\mu \tau}\,,
$$
for $\tau\ge\tau_1$. As a consequence,  
$$
\|u(t)\|_{L^1} \,=\, \|\hat v(\ln(1{+}t)-\tau_0)\|_{L^1} \,\le\,  
C_{10}\,e^{\mu\tau_0}\,(1+t)^{-\mu}\quad \mbox{for } t \ge
e^{\tau_1+\tau_0}-1\,.  
$$
By (\ref{b2}), we also have $\|u(t)\|_{L^\infty} \le
C\,t^{-\mu-(N/2)}$ for $t$ sufficiently large, a property which
implies that $u\equiv 0$ by \cite[Proposition~3]{BLS02} and 
\cite[Corollary~4.2]{GGK03}. This contradicts the fact that 
$u(t,x) > 0$ for $t > 0$, hence we have proved the first assertion 
in (\ref{d3}). As for the second claim, it is a straightforward 
consequence of (\ref{b0}) and (\ref{d4}). 

Now, since $\|\hat v(\tau)\|_m \to 0$ as $\tau \to \infty$, 
it follows from (\ref{d2}) that there exists $\tau_2 \ge 0$ such that 
$\|w(\tau)\|_m \le \varrho$ for all $\tau\ge\tau_2$. On the other
hand, as $w(\tau)\in W_c$ for each $\tau\ge 0$, we have 
$w(\tau,\xi) = M(\tau)\,G(\xi) + f(M(\tau)\,G(\xi))$ for $(\tau,\xi)
\in [0,\infty)\times\RR^N$, where $f$ is as in Theorem~\ref{thc7}. 
In view of (\ref{c2}) and (\ref{c4}) we deduce that, 
for $\tau\ge\tau_2$, 
\beqn
\label{d5}
\frac{dM}{d\tau}(\tau) \,=\, -\int_{\RR^N} |\nabla w(\tau,
\xi)|^{q_\star}\dd\xi \,=\,  -\|\nabla G\|_{L^{q_\star}}^{q_\star}
\,M(\tau)^{q_\star} - \omega(\tau)\,, 
\eeqn
where
$$
\omega(\tau) \,=\, \int_{\RR^N} \Bigl(|\nabla
w(\tau,\xi)|^{q_\star} - M(\tau)^{q_\star}\, |\nabla
G(\xi)|^{q_\star}\Bigr)\dd\xi\,. 
$$
To bound $\omega(\tau)$, we remark that $|\,|y+z|^{q_\star}-
|y|^{q_\star}\,| \le {q_\star}(|y|+|z|)^{q_\star-1}|z|$ for 
all $y,z \in \RR$. Also, since 
$$
\frac{1}{2} + \frac{q_\star-1}{2} + \frac{N}{2(N+1)} \,=\, 1\,, \quad 
\hbox{and}\quad 2mq_\star\,\frac{N+1}N \,=\, 2m\,\frac{N+2}N \,>\, 
N+2\,,
$$
it follows from H\"older's inequality that, for all $g,h \in 
L^2_m(\RR^N)$, 
$$
  \|\,|h|^{q_\star-1} g\|_{L^1} \,\le\, \|b^m h\|_{L^2}^{q_\star-1}
  \,\|b^m g\|_{L^2} \,\|b^{-mq_\star}\|_{L^{2(N+1)/N}} \,\le\, 
  C \,|h|_m^{q_\star-1} |g|_m\,,
$$
where $b(\xi) = (1+|\xi|^2)^{1/2}$. Thus 
\bean
|\omega(\tau)|& \le & q_\star\,\int_{\RR^N} \Bigl(M(\tau)\,
|\nabla G(\xi)| + |\nabla f(M(\tau)\,G(\xi))|\Bigr)^{q_\star-1}\, 
|\nabla f(M(\tau)\,G(\xi))|\dd\xi \\
& \le & C\,(M(\tau)\,\|G\|_m + \|f(M(\tau)\, G)\|_m)^{q_\star-1}
  \,\|f(M(\tau)\, G)\|_m \phantom{\frac12}\\ 
& \le & C\,M(\tau)^{q_\star-1}\, \|f( M(\tau)\, G)\|_m\,, 
\eean
where in the last inequality we have used the fact that $f(0)= 0$ 
and $f$ is globally Lipschitz continuous. Since
$f\in\mathcal{C}^1(E_c;E_s)$ with $f(0)=0$ and $Df(0)=0$, the above
inequality and (\ref{d3}) imply that 
\beqn
\label{d6}
\lim_{\tau\to\infty} \frac{\omega(\tau)}{M(\tau)^{q_\star}} \,=\, 0\,.
\eeqn
Combining (\ref{d3}), (\ref{d5}), (\ref{d6}), and recalling that
$q_\star-1 = 1/(N{+}1)$, we conclude that
\beqn
\label{d7}
\lim_{\tau\to\infty} \tau M(\tau)^{q_\star-1} \,=\, \frac{1}
{(q_\star{-}1) \|\nabla G\|_{L^{q_\star}}^{q_\star}}\,, 
\quad \hbox{i.e.} \quad 
\lim_{\tau\to\infty} \tau^{N+1} M(\tau) \,=\, M_\star\,,
\eeqn
where $M_\star$ is as in Theorem~\ref{th1}. As $w(\tau,\xi) =
M(\tau)\,G(\xi) + f(M(\tau)\,G(\xi))$, we deduce from (\ref{d3}), 
(\ref{d7}) and the properties of $f$ that $\tau^{N+1}\,\|w(\tau) - 
M(\tau)\,G\|_m\to 0$ as $\tau\to\infty$. Combining this result with 
(\ref{d2}), (\ref{d7}), we arrive at  
\beqn
\label{d8}
\lim_{\tau\to\infty} \tau^{N+1}\,\left\|\hat v(\tau) -
\frac{M_\star}{\tau^{N+1}}\, G \right\|_{L^1} \,=\, 0\,. 
\eeqn
Of course, the same result holds for $v(\tau) = \hat v(\tau-\tau_0)$.
If we now return to the original function $u(t,x)$ via the 
transformation (\ref{newuf}), we obtain exactly (\ref{dr}) 
for $p=1$. The case $p \in (1,\infty]$ then follows from (\ref{b2}) 
by a classical interpolation argument. \qed

\end{document}